\documentclass[default]{sn-jnl}
\usepackage{amsfonts}
\usepackage{graphicx}

\jyear{2021}
\theoremstyle{thmstyleone}

\theoremstyle{thmstyletwo}
\newtheorem*{result}{{\bf Main Result}}
\theoremstyle{thmstylethree}

\def\r{\mathbb R}
\def\e{\mathbb E}

\begin{document}
\title{Effect of local fractional derivatives on Riemann curvature tensor}
 
\author{\fnm{Muhittin Evren} \sur{Aydin}}

\affil{\orgdiv{Department of Mathematics, Faculty of Science}, \orgname{Firat University},   \city{Elazig}, \postcode{23200},   \country{Turkey}, meaydin@firat.edu.tr}

\abstract{In this short note, we investigate the effect of the local fractional derivatives on the Riemann curvature tensor that is a common tool in calculating curvature of a Riemannian manifold. For this, first we introduce a general local fractional derivative operator that involves the mostly used ones in the literature as conformable, alternative, truncated $M-$ and $\mathcal{V}-$fractional derivatives. Then, according to this general operator, a particular Riemannian metric tensor field on the real affine space $\r^{n}$ that is different than Euclidean one is defined. In conlusion, we obtain that the Riemann curvature tensor of $\r^{n}$ endowed with this particular metric is identically $0$, namely, locally isometric to Euclidean space.}

 
\maketitle


\section{Introduction} \label{sec1}


Fractional calculus is based on the notions of differentiation and
integration of real or complex orders and laid its foundations at near time when the usual calculus emerges. Since the middle of 17th
century, the first class mathematicians of the history, as Leibniz,
L'Hopital, Euler, Bernoulli, Riemann and among others, made a contribution
to this field. We refer to \cite{bf}, \cite{skm} for the structural and historical developments.

Besides much interest in the pure mathematics is given to fractional calculus, applications appear in the natural sciences when a process or system need to be characterized with spatial and temporal nonlocality. For example, those applications rang from medicine (\cite{dcbb}, \cite{nausaa}), bioengineering (\cite{m}), viscoelasticity (\cite{yn}) to dynamical systems (\cite{yn1}).

Contrary to the derivative and integral of integer order, many different
versions of the fractional derivative and integral exist. While
the mostly used non-local fractional derivatives are Riemann-Liouville, Caputo, Riesz fractional derivatives among others (\cite{p}, \cite{skm}), the local
ones are conformable (\cite{a}, \cite{kays}), alternative (\cite{k}), truncated $M-$ and $\mathcal{V}-$fractional derivatives (\cite{v}, \cite{v1}).

In this paper, our purpose is to investigate the effect of a fractional derivative operator on curvature of the real affine $n-$space $\r^{n}$. Explicitly, we concentrate on the Riemann curvature tensor $R$ playing a crucial role in Riemannian geometry because it is a common tool to calculate curvature of a Riemannian manifold. For example, for a Riemannian manifold $M$ if $R$ is identically $0$ then $M$ is locally isometric to $\r^n$, or equivalently, each point of $M$ has a neighborhood that is isometric to an open subset in $\r^n$ with Euclidean metric (see \cite[Theorem 7.3]{l}). In this sense, $M$ is also said to be {\it flat}. In addition, when the sectional curvature is identically a positive (resp. negative) constant then $M$ is locally isometric to a sphere (resp. hyperbolic space), which those are known as {\it model spaces} of Riemannian geometry (\cite{on}).

In the case of non-local fractional derivatives, the main features that the usual derivative of integer order has are not satisfied, for example the usual Leibniz and chain rules (\cite{skm}, \cite{t}, \cite{t1}, \cite{t2}). The absence of those features reveals
a big difficulty if one needs to set a theory in differential geometry
because of their essential roles. In other words, if one considers a non-local fractional derivative instead of the usual derivative of integer order,
one would not be able to use the techniques from Riemannian geometry due to
the absence of Leibniz and chain rules. Numerous notions in Riemannian
geometry are based on those rules.

In contrast to the non-local fractional derivatives, the local fractional derivatives have the same main features that the usual derivative
of integer order holds as Leibniz rule, chain rule, mean value theorem,
Rolle Theorem, and so on (see \cite{v1}). This is because it is reasonable to
tend to a local fractional derivative when a differential geometrical
problem from fractional calculus viewpoint is considered.

The organization of this paper is as follows. In Sect. \ref{sec2} we introduce a general local fractional derivative operator (see Eq. \eqref{intro3}) involving conformable, alternative, truncated $M-$ and $\mathcal{V}-$fractional derivatives. In Sect. \ref{sec3}, we recall some basic notions from the Riemannian geometry.  In Sect. \ref{sec4}, we present a main result \ref{ex} which points out that the addressed local fractional derivatives do not affect the Riemann curvature tensor. In order to do this, a certain Riemannian metric $\langle , \rangle_{\alpha}$ on $\r^{n}$ that is different than Euclidean one is introduced so that this general operator is an orthonormal coordinate vector field with respect to $\langle , \rangle_{\alpha}$ ($0 < \alpha \leq 1$). Then we obtain that $(\r^{n},\langle , \rangle_{\alpha})$ is again flat independently from the value of $\alpha$. 

\section{Local fractional derivatives} \label{sec2}

We recall the definition of $\mathcal{V}-$ fractional derivative introduced in \cite{v1} because it is a general version of the local fractional derivatives we addressed in Section \ref{sec1}. See also the nice diagram on page 23 of the cited paper. Hence when a statement holds for $\mathcal{V}-$ fractional derivative then so directly do for the others.

Denote by $\r_+$ the set of all the positive real numbers. Let $\gamma , \beta , \rho , \delta $ be some complex numbers with positive real parts and $p,q \in \r_+$. Introduce a {\it six-parameter Mittag–Leffler function} (\cite{v1})
$$ _i  \e_{\gamma , \beta , p}^{\rho , \delta , q}(z)=  \sum_{k=0}^{i}\frac{(\rho)_{qk}}{(\delta)_{pk}} \frac{z^k}{\Gamma(\gamma k +\beta)}, \quad z \in \mathbb{C}, Re(z)>0,$$
where $\Gamma (\rho )$ is the gamma function on $\rho$ and $(\rho)_{qk}$ is a generalization of the Pochhammer symbol defined by $(\rho)_{qk}=\Gamma(\rho +q k)/\Gamma(\rho )$.

Let $0 < \alpha \leq 1$ be a real number and $t \mapsto f(t) $, $t \in I \subset \r_+$. Then the truncated $\mathcal{V}-$fractional derivative of $f(t)$ is
\begin{equation}
 _i^\rho \mathcal{V}_{\gamma , \beta , \alpha}^{\delta , p , q} f(t) = \lim_{\epsilon \to 0} \frac{f(t _i H_{\gamma , \beta , p}^{\rho , \delta , q}(\epsilon t^{-\alpha})) -f(t)}{\epsilon}, \label{intro1}
\end{equation}
where $ _i H_{\gamma , \beta , p}^{\rho , \delta , q}(\epsilon t^{-\alpha}) (z)= \Gamma (\beta) _i  \e_{\gamma , \beta , p}^{\rho , \delta , q}(z) $. 

A function is said to be $\alpha-${\it differentiable} if the limit in Eq. \eqref{intro1} exists. In addition, if $f(t)$ is differentiable (in the usual sense) in a neighborhood of some $a>0$ and if $\lim_{t \to 0^+} (_i^\rho \mathcal{V}_{\gamma , \beta , \alpha}^{\delta , p , q} f(t))$ exists, then 
$$ _i^\rho \mathcal{V}_{\gamma , \beta , \alpha}^{\delta , p , q} f(0) = \lim_{t \to 0^{+}} (_i^\rho \mathcal{V}_{\gamma , \beta , \alpha}^{\delta , p , q} f(t)).$$

The truncated $\mathcal{V}-$fractional derivative is a linear operator and satisfies Leibniz and chain rules as well as (\cite{v1})
\begin{equation}
  _i^\rho \mathcal{V}_{\gamma , \beta , \alpha}^{\delta , p , q} f(t) = \frac{t^{1-\alpha}\Gamma(\beta)(\rho)_q}{\Gamma(\gamma + \beta)(\delta)_p} \frac{df(t)}{dt}. \label{intro2}
\end{equation}
According to the particular values of the parameters $\gamma , \beta , \delta , p , q$, the other local fractional derivatives as conformable,  alternative, truncated $M-$fractional derivatives and etc., can be derived.

As can be deduced from Eq. \eqref{intro2}, the truncated $\mathcal{V}-$fractional derivative is linearly related to the usual derivative of integer order. Hence, a general local fractional derivative operator is defined as follows.

Let a positive real-valued function $c_{\alpha }(t)$ of class $C^{2}$ be given by
\begin{equation*}
t\mapsto c_{\alpha }(t)\in \mathbb{R}_{+},\quad t\in I\subset \mathbb{R}_{+}.
\end{equation*}
Since the coefficient of the term $df(t)/dt$ in \eqref{intro2} is function of $t$, we now introduce generally%
\begin{equation}
\frac{d^{\alpha }}{dt^{\alpha }}:=c_{\alpha }(t)\frac{d}{dt}, \label{intro3}
\end{equation}%
where as usual $d/dt$ is the derivative operator of first order. For instance, Eq. \eqref{intro3} implies conformable fractional derivative if $c_{\alpha }(t)=t^{1-\alpha }$ and a truncated $M-$fractional derivative if $%
c_{\alpha }(t)=t^{1-\alpha }/\Gamma(1+\beta)$ ($\beta \in \r_+$) and etc.

\section{Curvature tensor field} \label{sec3}

Let $(M, \langle , \rangle)$ be a $n-$dimensional Riemannian manifold and $\nabla$ the Levi-Civita connection on $M$. Denote by $(x_1 , ... , x_n)$ a coordinate system on a neighborhood in $M$ such that $\partial_i =\partial / \partial x_i$ $(1\leq i \leq n)$. 

The {\it Christoffel symbols} $\Gamma_{ij}^{k}$ for the system $(x_1 , ... , x_n)$ are (\cite{on})
$$ \nabla_{\partial_i}\partial_j =\sum_{k=1}^n \Gamma_{ij}^{k} \partial_k, \quad 1\leq i,j \leq n.$$
If we denote by $g_{ij}$ the components of the Riemannian metric $\langle , \rangle$, then
\begin{equation}
\Gamma_{ij}^{k} = \frac{1}{2}  \sum_{m=1}^n g^{km} \left ( \frac{\partial g_{jm}}{\partial x_i} +\frac{\partial g_{im}}{\partial x_j} -\frac{\partial g_{ij}}{\partial x_m} \right ), \quad 1\leq i, j,k \leq n , \label{tensor1}
\end{equation}
where $(g^{ij})$ is the inverse matrix of $(g_{ij})$.

The {\it Riemann curvature tensor} $R$ of $(M , \left\langle ,\right\rangle)$ is
$$R(X,Y)Z = \nabla_{[X,Y]}Z - \nabla_{X}\nabla_{Y}Z + \nabla_{Y}\nabla_{X}Z ,$$
where $X,Y,Z$ are arbitrary smooth vector fields on $M$. Denoting $R_{jkl}^i$ the components of $R$ $(1\leq i,j,k,l \leq n)$, 
\begin{equation*}
 R(\partial_k,\partial_l)\partial_j  = \sum_{i=1}^n R_{jkl}^i \partial_i. \label{tensor2}
\end{equation*}
In terms of the Christoffel symbols,
\begin{equation}
R_{jkl}^i = \frac{\partial }{\partial x_l} \Gamma_{kj}^i - \frac{\partial }{\partial x_k} \Gamma_{lj}^i
+ \sum_{m=1}^n \Gamma_{lm}^i \Gamma_{kj}^m - \sum_{m=1}^n \Gamma_{km}^i \Gamma_{lj}^m ,\quad 1\leq i,j,k,l \leq n. \label{tensor2}
\end{equation}
Notice here that, in the Euclidean case, $R_{jkl}^i$ are identically $0$. 

\section{Main result} \label{sec4}

Let $(\mathbb{R}^{n},\left\langle ,\right\rangle _{e})$\ be Euclidean $n-$%
space and $(x_{1},...,x_{n})$ the usual coordinate system. Denote by $%
\left\vert \cdot \right\vert _{e}$ the Euclidean norm induced by the
Euclidean metric $\left\langle ,\right\rangle _{e}.$ Let $\partial_i$ and $dx_{i}$ be the $i-$th coordinat vector field and its dual, for $1\leq i \leq n.$ Then, $\left\langle ,\right\rangle
_{e}=dx_{1}^{2}+...+dx_{n}^{2}$, where $dx_i^2=dx_i \otimes dx_i$ for the tensor product $\otimes$.

Assume that $x_{i}\in I_{i}\subset \mathbb{R%
}_{+}$ ($1 \leq i \leq n$). The multidimensional version of Eq. \eqref{intro3} is 
\begin{equation}
\frac{\partial ^{\alpha }}{\partial x_{i}^{\alpha }}:=c_{\alpha }(x_{i})%
\frac{\partial }{\partial x_{i}},  \label{intro4}
\end{equation}%
where $c_{\alpha }(x_{i})$ is a positive single-valued function of class $C^{2}$ on $I_{i}.$ 

Let $f(x_1,...,x_n)$ be a real-valued function on $I:=I_1 \times ... \times I_n$. If $\partial ^{\alpha }f/\partial x_{i}^{\alpha }$ exists on $I$ for every $1 \leq i \leq n$ then it is said to be $\alpha-${\it differentiable} on $I$. Since $c_{\alpha }(x_{i})$ is of class $C^{2}$ on $I_{i}$, the $\alpha-$differentiability of $f(x_1,...,x_n)$ on $I$ is equivalent to be of class $C^1$ on $I$.

Let $ \partial_i^{\alpha} = \partial ^{\alpha }/\partial x_{i}^{\alpha }$. Because $\left\vert \partial_i^{\alpha}\right\vert_{e}=c_{\alpha }(x_{i}),$ $\partial_i^{\alpha }$ is
not unitary for any $1\leq i \leq n$, and so we need to discuss a new metric tensor field
besides the Euclidean one. For this, we introduce%
\begin{equation}
\left\langle ,\right\rangle _{\alpha }:=\frac{dx_{1}^{2}}{c_{\alpha
}(x_{1})^{2}}+...+\frac{dx_{n}^{2}}{c_{\alpha }(x_{n})^{2}},  \label{intro5}
\end{equation}%
which is obviously a Riemannian
metric because $c_{\alpha }(x_{i})\in \mathbb{R}_{+}$, for every $1\leq i \leq n$. Hence, we may deduce $\langle \partial_i ^{\alpha } , \partial_j^{\alpha } \rangle_{\alpha } =\delta_{ij}$, $ 1 \leq i,j \leq n,$ where $\delta_{ij}$ is the Kronecker's delta.

Since $\mathbb{R}^{n}$ is a differentiable manifold, $(\mathbb{R}%
^{n},\left\langle ,\right\rangle _{\alpha })$ is now a Riemannian manifold and
hence we are able to calculate the Riemann curvature tensor of $(%
\mathbb{R}^{n},\left\langle ,\right\rangle _{\alpha })$, denoted by $R^{\alpha}$.  The following result points out the ineffectiveness of $\partial_i^{\alpha}$ on $R^{\alpha}$.
\begin{result}\label{ex}
The Riemannian manifold $(\r^n,\left\langle ,\right\rangle _{\alpha })$ is flat, independently from the value of $0 < \alpha \leq 1$, where $\left\langle ,\right\rangle _{\alpha }$ is defined by Eq. \eqref{intro5}.
\end{result}

\begin{proof}[Proof of Main Result]
We will derive the tensor field $R^{\alpha}$ of $(\r^n,\left\langle ,\right\rangle _{\alpha })$. This is equivalent to compute the components $R_{jkl}^{\alpha ,i}$, where $0 < \alpha \leq 1$. For this, first we will find the Christoffel symbols $\Gamma_{ij}^{\alpha, k}$ of $(\r^n,\left\langle ,\right\rangle _{\alpha })$. The components of $\left\langle ,\right\rangle _{\alpha }$ are 
$$g_{ij}^{\alpha}=c_{\alpha}(x_{i})^{-2}\delta_{ij}, \quad g^{\alpha , ij}=c_{\alpha}(x_{i})^{2}\delta_{ij}, \quad 1\leq i,j \leq n.$$
Therefore, because $g^{\alpha , ij}=0$ for $1\leq i \neq j \leq n$, Eq. \eqref{tensor1} is now
\begin{equation}
\Gamma_{ij}^{\alpha, k} = \frac{1}{2}   g^{\alpha, kk} \left ( \frac{\partial g_{jk}^{\alpha}}{\partial x_i} +\frac{\partial g_{ik}^{\alpha}}{\partial x_j} -\frac{\partial g_{ij}^{\alpha}}{\partial x_k} \right ), \quad 1\leq i, j,k \leq n . \label{tensor3}
\end{equation}
The parenthesis in Eq. \eqref{tensor3} is nonzero only if $i=j=k$. For any other case, $\Gamma_{ij}^{\alpha, k}$ are identically $0$. Consequently, the nonzero Christoffel symbols are
$$ \Gamma_{ii}^{\alpha, i} = \frac{1}{2}   g^{\alpha, ii}  \frac{\partial g_{ii}^{\alpha}}{\partial x_i}  , \quad 1\leq i \leq n , $$
or equivalently,
\begin{equation}
 \Gamma_{ii}^{\alpha, i} =  - \frac{1}{ c_{\alpha}(x_{i})} \frac{dc_{\alpha}}{dx_i}(x_{i}) , \quad 1\leq i \leq n . \label{tensor4} 
\end{equation}

The components $R_{jkl}^{\alpha, i}$ are skew-symmetric up to the indices $j,k$, which it must be $j\neq k$. Hence, because Eqs. \eqref{tensor2} and \eqref{tensor4}, we have
$$ R_{jkl}^{\alpha, i} = - \frac{\partial }{\partial x_k} \Gamma_{lj}^{\alpha, i}
 - \sum_{m=1}^n \Gamma_{km}^{\alpha, i} \Gamma_{lj}^{\alpha, m}  , \quad 1\leq i, j \neq k , l \leq n ,$$
where the partial derivative must be zero due to $ j \neq k$. Analogously, all the terms of the sum are again zero, obtaining $R_{jkl}^{\alpha, i}=0$, for $1\leq i, j \neq k , l\leq n$.
\end{proof}

\section{Conclusions}

In this note, by concentrating on the Riemann curvature tensor we investigated how the local fractional derivatives defined in the literature affect curvature of $\r^n$. This is because of the crucial role of the curvature tensor field that we addressed in Section \ref{sec1}.

As general, we considered the local fractional derivative operator $\partial ^{\alpha }/\partial x_{i}^{\alpha }$ given by Eq. \eqref{intro4} ($0 < \alpha \leq 1$), which it includes the well-known local fractional derivatives, i.e. conformable, alternative, truncated $M-$ and $\mathcal{V}-$fractional derivatives. According to $\partial ^{\alpha }/\partial x_{i}^{\alpha }$, we introduced a different Riemannian metric tensor field $\left\langle ,\right\rangle _{\alpha }$ (see Eq. \eqref{intro5}) on $\r^n$ than Euclidean one such that $\langle \partial ^{\alpha }/\partial x_{i}^{\alpha } , \partial^{\alpha }/\partial x_{j}^{\alpha } \rangle_{\alpha } =\delta_{ij}$ ($1 \leq i,j \leq n$). 

Our main result \ref{ex} states that the Riemann curvature tensor $R^{\alpha}$ of $(\r^n , \left\langle ,\right\rangle _{\alpha })$ is identically $0$, that is, $(\r^n , \left\langle ,\right\rangle _{\alpha })$ is locally isometric to Euclidean space, independently from the value of $\alpha$. From this result, we understood that none of the addressed local fractional derivatives has an effect on the geometry of $\r^n$.

\end{document}